\documentclass{article}

\usepackage{amsmath}
\usepackage{amsfonts}
\usepackage{amssymb}
\usepackage{mathrsfs}
\usepackage{graphicx}
\usepackage{xcolor}

\newcommand   \primeInt {\mathtt{q}}
\newcommand   \fieldChar {\mathtt{\,p\,}}
\newcommand   \finiteIntegerField {\basel{\mathbb Z}{\fieldChar}\,} 
\newcommand   \scalars  {\mathbb F}

\renewcommand   \mod  {~\textsf{mod}~ }
\renewcommand   \gcd  {~\textsf{gcd}~ }

\renewcommand   \det  {\textsf{det}}

\newcommand  \xx   {\mathbf x}

\newcommand  \polynomials[3] { #1{\mathbf{[}}\basel{#2}{\mathrm{1}},\, \ldots,\, \basel{#2}{\mathit{#3}}{\mathbf{]}} }
\newcommand  \singlevariablepolynomials[2] {#1{\mathbf{[}#2\mathbf{]}}}
\newcommand  \basel[2]{#1_{_{#2}}}

\newcommand  \tab  {\hspace*{0.5cm}}
\newcommand  \ltab  {\hspace*{-0.5cm}}

\newcommand  \shiftright  {\hspace*{2.0cm}}
\newcommand  \shiftleft  {\hspace*{-2.0cm}}

\newcommand  \bglb {\big (}
\newcommand  \bgrb {\big )}

\newcommand    \NSPACE  {\mathcal {NPSPACE}}

\newcommand    \PSPACE  {\mathcal {PSPACE}}

\newcommand    \IP  {\mathcal {IP}}

\newcommand    \BoundedErrorProbPTIME  {\mathcal {BPP}}

\newcommand    \PTIME  {\mathcal {P}}

\newcommand    \NPTIME  {\mathcal {NP}}

\newcommand    \NSPACEwithProof  {\mathcal{NPSPACE}\mathsf{\_with\_proof\_in\_}\PSPACE}

\newcommand    \CoNPTIME  {\textsf{Co-}\mathcal {NP}}

\definecolor{titlecolor}{RGB}{144,48,48}
\definecolor{authorname}{RGB}{16,96,16}%
\definecolor{addrscolor}{RGB}{60,113,183}
\definecolor{secheader}{RGB}{16,32,128}%
\definecolor{refscolor}{RGB}{16,64,128}%

\begin{document}

\title{\textcolor{titlecolor}{\bf{Complexity of Solution of Simultaneous Multivariate Polynomial Equations}}}

\author{\\
\\
{\textcolor{authorname}{\bf{Duggirala Meher Krishna}}}\\
{\textcolor{addrscolor}{\small{\bf{Gayatri Vidya Parishad College of Engineering (Autonomous)}}}} \\
{\textcolor{addrscolor}{\small{\bf{Madhurawada, VISAKHAPATNAM -- 530 048, Andhra Pradesh, India}}}} \\
{\textcolor{addrscolor}{\small{\bf{E-mail ~: \tab duggiralameherkrishna@gmail.com}}}}\\
 \\
 {\textcolor{addrscolor}{and}} \\
 \\
{\textcolor{authorname}{\bf{Duggirala Ravi}}}\\
{\textcolor{addrscolor}{\small{\bf{Gayatri Vidya Parishad College of Engineering (Autonomous)}}}} \\
{\textcolor{addrscolor}{\small{\bf{Madhurawada, VISAKHAPATNAM -- 530 048, Andhra Pradesh, India}}}} \\
\shiftleft \tab {\textcolor{addrscolor}{\small{\bf{E-mail ~: \tab ravi@gvpce.ac.in; \tab duggirala.ravi@yahoo.com};}}} \\
{\textcolor{addrscolor}{\small{\shiftright \ltab \ltab \bf{ duggirala.ravi@rediffmail.com; \tab drdravi2000@yahoo.com}}}}
\\
\\
}

\date{}

\maketitle

\textcolor{secheader}{
\begin{abstract}
In this paper, an original reduction algorithm for solving simultaneous multivariate polynomial equations is presented. The algorithm is exponential in complexity, but the well-known algorithms, such as the extended Euclidean algorithm and Buchberger's algorithm, are superexponential. The superexponential complexity of the
well-known algorithms is due to their not being ``minimal'' in a certain sense. Buchberger's algorithm produces
a Gr\"{o}bner basis. The proposed original reduction algorithm achieves the required task {\em via} computation
of determinants of parametric Sylvester matrices, and produces a Rabin basis, which is shown to be minimal, when two multivariate polynomials are reduced at a time. The minimality of Rabin basis allows us to prove exponential lower bounds for the space complexity of an algebraic proof of certification, for a specific computational problem in the computational complexity class $\PSPACE$, showing that the complexity classes $\PSPACE$ and $\PTIME$ cannot be the same. 
By the same reasoning, it follows that  $\CoNPTIME \neq  \NPTIME$ and $\CoNPTIME \neq \PTIME$, and that
the polynomial time hierarchy does not collapse. It is also shown that the class $\BoundedErrorProbPTIME$
of languages decidable by bounded error probabilistic algorithms with (probabilistic) polynomial time proofs
for the membership of input words is not the same as any one of the complexity classes $\PTIME$, $\NPTIME$
and $\CoNPTIME$. By relativization of $\BoundedErrorProbPTIME $ with respect to $\PTIME$ and $\NPTIME$,
it can be again inferred that $\NPTIME \neq \PTIME$ and that $\CoNPTIME \neq \PTIME$.
\end{abstract}
}

\begin{small}
\textcolor{secheader}{
\paragraph{Keywords.}~~Extended Euclidean algorithm; ~Buchberger's algorithm; ~Gr\"{o}bner basis; ~ Parametric Sylvester matrices; ~ Parametric resolvant; ~ Rabin basis; ~ Minimality; ~ Algebraic complexity theory.
}
\end{small}
\\


{\textcolor{titlecolor}{\textbf{\section{Introduction}}}}
Two prominent methods for reducing simultaneous multivariate polynomial equations are the extended
Euclidean algorithm and Buchberger's algorithm. The reduction can be performed eliminating one variable
at a time. If two multivariate polynomials vanish simultaneously, then so does their parametric $\gcd$,
where the parametric$\gcd$ is an element in the integral domain or the field of rational functions in
the remaining variables. However, the converse may not be true: for two multivariate polynomials,
if the parametric $\gcd$ vanishes, for a particular interpretation of variables from the
algebraic closure of the ground field without any free variables, then the same cannot be
guaranteed for the two multivariate polynomials, whose $\gcd$ has vanished under the chosen
interpretation. The failure of the converse to hold true in general contributes to the
superexponential complexity of these two well-known algorithms (Section \ref{Sec-EEA-BA-GB}).

Now, for the two multivariate polynomials in discussion to vanish simultaneously, it is
both necessary and sufficient that the determinant of the parametric Sylvester matrix,
called their {\em {resultant}} with respect to the variable being eliminated, vanishes,
for any or some values in the ground field, in which the zeros are being searched for,
assuming that neither of the two multivariate polynomials vanishes identically, under
the chosen interpretation nullifying the resultant. The equivalence of the vanishing
resultant --- except when at least one of the two multivariate polynomials identically
vanishes, under the interpretation of values to the variables other than the variable being
eliminated in the current reduction step --- to the sharing of a common zero in the algebraic
closure of the ground field without parameters lends us the minimality criterion, for the
reduction of the two multivariate polynomials. The entries of the Sylvester matrix being
multivariate polynomials themselves, {\em {albeit}} without the variable being eliminated
in the current step, the Gaussian elimination procedure cannot be applied, even though the
final resultant, which is the determinant of the Sylvester matrix, is the same. The reason
for inapplicability of the Gaussian elimination procedure for the determinant of the 
parametric matrices is that the computations in the intermediate stages might become
excessively large, causing superexponential space and time for the completion of its
computations, as can be achieved by simplification of intermediate results. The overall
performance of the Gaussian elimination procedure for computation of the parametric
resultant may be comparable to the extended Euclidean algorithm or Buchberger's algorithm,
if not worse than either. A step-by-step simplified computation of the resultant that
does not run into the space or time explosion problem, which is experienced with
the Gaussian elimination method for the parametric matrices, is also presented.
The reduced multivariate polynomial basis obtained by taking the resultant for
each reduction step is called a {\em {Rabin basis}}, in honour of
Professor Michael Oser Rabin, for his profound contributions to computer science
(Section \ref{Sec-SM-Res-RB}). 

The minimality criterion allows us to derive exponential lower bounds for space requirement
for an algebraic proof of certification, for a specific computational problem in two variables
with undetermined coefficients over any finite field. The specific computational problem is
shown to exist in $\PSPACE$, by exhibiting an algorithm for solving it, requiring space bounded
by a linear expression in the sum of the degrees of the two independent variables and logarithm
of the cardinality of the finite field, excluding the space required for the finite field
arithmetic operations.  It is customarily acknowledged that, for a computational problem
to be in $\PTIME$, a polynomial time algorithm, together with a proof of certification
--- the proof being bounded in space by a polynomial in appropriate parameter values for
the instance --- must exist for its correctness of operation. The class of nondeterministic
computational or decision problems, for which it is possible to produce machine checkable
deterministic proofs, {\em
{bounded in space by a polynomial for each such specific
computational problem}}, denoted by $\NSPACEwithProof$, relative to any particular fixed
system of deductive or symbolic logic, equipped with rules of inference, that might be
extensible, is included in $\NPTIME$, and hence the two complexity classes $\NSPACEwithProof$
$\NPTIME$ represent the same complexity class.
The machine checkable proofs may include references to external facts, the rules of
inference may be specialized to a specific computational or decision problem, and the
extensibility is the system's or users' ability to add more rules perhaps adaptively 
and / or interactively. The nonexistence of a polynomial time deterministic verification
algorithm for a computational problem can be inferred from the nonexistence of a proof
of correctness for any such algorithm, for its solution, that is bounded in space by
a polynomial in the acceptable parameter values for its instances. A discernment of
Herbrand's theorem, as applied to multivariate polynomials, shows that there cannot
be a shorter form for the algebraic proof of certification, for the specific
computational problem in $\PSPACE$, because an immediate reflection shows that
the degree of the resultant for an instance to the computational problem under
investigation, even when the number variables is only two, is exponential, which
must be combined with the fact that further $\mod$ and $\gcd$ operations may have
to be performed, for the completion of the proof. In particular, the algebraic form
of the resultant needs to hold for all undetermined coefficients, degrees and
field characteristics, for the application of Herbrand's theorem. Moreover, it is
easily possible to assume multivariate polynomials with undetermined coefficients
(with more than just two independent variables and undetermined coefficients) as instances
for the specific computational problem in $\PSPACE$, and again invoking Herbrand's theorem,
recursion can be applied, to produce an algebraic proof of correctness for the specific
computational problem in discussion. By restricting interpretation of variables to small
dimension extension fields, a deterministic algorithm running in linear space (possibly
excluding the space required for the finite field arithmetic operations) can be exhibited,
for the specific computational problem with generalization to multiple variables. The
occurrence of recursion effectively annihilates any little hope of finding a proof of
certification bounded in space by a polynomial, for any deterministic algorithm for the
specific computational problem in its most generality, even when the interpretation is
restricted to small dimension extension fields. In summary, we have to become contented
in accepting that $\PSPACE$ cannot be $\PTIME$. In fact, an almighty can be assumed
to be capable of guessing the correct answer to the question posed as part of the
specific computational problem, but the impossibility of producing a polynomial
space proof of certification shows that $\PSPACE$ cannot be $\NPTIME$, either.
By the same reasoning, it follows that $\CoNPTIME \neq \PTIME$ and
$\CoNPTIME \neq \NPTIME$, and that the polynomial time hierarchy does not collapse.
In addition to these results, it is also shown that the class  $\BoundedErrorProbPTIME$
of languages acceptable by bounded error probabilistic algorithms with probabilistic
polynomial time proofs for the membership of an input word is not the same as any one of
the complexity classes $\PTIME$, $\NPTIME$ and $\CoNPTIME$.  From the discussions,
it follows again that $\NPTIME \neq \PTIME$ and that $\CoNPTIME \neq \PTIME$,
by relativization of $\BoundedErrorProbPTIME $, with respect to $\PTIME$ and $\NPTIME$
(Section \ref{Sec-PSPACE-not-equal-to-NP-or-P}).

{\textcolor{titlecolor}{\textbf{\section{\label{Sec-EEA-BA-GB}Extended Euclidean Algorithm, Buchberger's Algorithm and Gr\"{o}bner Basis}}}}

Let $\scalars$ be a field and $\polynomials{\scalars}{x}{n}$, for some positive integer $n \geq 2$,
be the integral domain of polynomials in $n$ independent variables $\basel{x}{1}, \, \ldots,\, \basel{x}{n}$,
with coefficients in $\scalars$. Let $\xx = (\basel{x}{1},\, \ldots, \, \basel{x}{n})$,
$\alpha(\xx) = \sum_{i = 0}^{d} \basel{a}{i}(\basel{x}{1}, \, \ldots,\, \basel{x}{n-1}) \basel{x^{i}}{n}$ and  $\beta(\xx) = \sum_{i = 0}^{d} \basel{b}{i}(\basel{x}{1}, \, \ldots,\, \basel{x}{n-1}) \basel{x^{i}}{n}$ be two polynomials in $\polynomials{\scalars}{x}{n}$, both of 
degree $d \geq 1$. It is further assumed that the polynomials $\basel{a}{i}(\basel{x}{1}, \, \ldots,\, \basel{x}{n-1}) $ and $\basel{b}{i}(\basel{x}{1}, \, \ldots,\, \basel{x}{n-1})$ in $\polynomials{\scalars}{x}{n-1}$
are all nonzero, and that each requires at least $\basel{L}{\min}$ units of space, for $0 \leq i \leq d$,
such that they could include more than $\basel{L}{\min} \geq 2$ terms with very diverse exponent vectors,
so that their products after expansion may contain only insignificantly small number of collision
terms, for applying cancellations or simplification of terms,  or they may admit succinct
representations requiring at least $\basel{L}{\min}$ units of space, when their products are not expanded.

The operation of the extended Euclidean algorithm for computation of the parametric $\gcd$ is
explained in the sequel. Since the two input polynomials are of the same degree $d$ in
$\basel{x}{n}$, an application of two consecutive steps to eliminate the highest degree term,
{\em {i.e.}}, $\basel{x^{d}}{n}$, results in two multivariate polynomials of degree $d-1$ each,
such that their coefficients would need $2\basel{L}{\min}$ units of space. Now, by induction,
an application of two consecutive steps to eliminate $\basel{x^{d-i}}{n}$, from the two
multivariate polynomials obtained as the result of the last consecutive pair of steps
by eliminating  $\basel{x^{d-i+1}}{n}$, for $i = 1,\, 2,\, ...\, d-1$, would result in
$2^{i}\basel{L}{\min}$ units of space, without expansion. Thus, when the products are not
expanded, the  overall space requirement for the elimination of
$\basel{x}{n}$ is at least ${\mathcal{O}}\bglb 2^{d} \basel{L}{\min} \bgrb$.
One more insight is concerning the final degree of any of the variables $\basel{x}{i}$,
for $1 \leq i \leq n-1$. For simplicity, let the degree of occurrence of the variable
$\basel{x}{i}$, for some fixed index $i$, where $1 \leq i \leq n-1$, be
$\basel{\delta}{i} \geq 2$, for each term occurring as the coefficient of
$\basel{x}{n}$ in either input polynomial. The elimination procedure produces
coefficients as multivariate polynomials in $\polynomials{\scalars}{x}{n-1}$.
Assuming that the occurrence of cancellations while simplifying the computations
is a rare event, the degree of occurrence of the variable $\basel{x}{i}$ in the
parametric $\gcd$ can be lower bounded by ${\mathcal O}\bglb 2^{d}\basel{\delta}{i}\bgrb$,
for $1 \leq i \leq n-1$.

Expansion and simplification of the products formed in the intermediate steps might not produce
a lot of cancellations, and would only be expected to further blow up the space requirement.
To eliminate $\basel{x^{d}}{n}$, by a consecutive pair steps, would need at least
${\mathcal{O}}\bglb \basel{L^{2}}{\min} \bgrb$ space, after expansion, and by induction,
to eliminate $\basel{x^{d-i+1}}{n}$, by a consecutive pair steps, would need at least
${\mathcal{O}}\bglb \basel{L^{2^{i}}}{\min} \bgrb$ space, after expansion,
for $i = 1, 2, ..., d$, resulting in the overall space requirement of at least
${\mathcal{O}}\bglb \basel{L^{2^{d}}}{\min} \bgrb$ space, after expansion. This is
the problem that causes the space explosion when the extended Euclidean algorithm is
applied, for eliminating the variable $\basel{x}{n}$, from the two input multivariate
polynomials $\alpha(\xx)$ and $\beta(\xx)$. Most of the zeros of the parametric $\gcd$
might not lead to the common zeros of $\alpha(\xx)$ and $\beta(\xx)$. 

Buchberger's algorithm follows closely the operational principle of the extended Euclidean algorithm
and, in effect, emulates the latter by considering the exponent vector as a whole, in the sum of terms
form. The multivariate polynomials so produced are collected in the Gr\"{o}bner basis \cite{Buchberger-1965},
named after the Ph D advisor of the author, presumably connoting Hilbert's basis theorem.
 
{\textcolor{titlecolor}{\textbf{\section{\label{Sec-SM-Res-RB}Parametric Sylvester Matrix, Parametric Resultant and Rabin Basis}}}}
Let $\alpha(\xx)$ and $\beta(\xx)$ be multivariate polynomials in $\polynomials{\scalars}{x}{n}$,
of degrees $\basel{d}{\alpha} \geq 1$ and $\basel{d}{\beta} \geq 1$. The Sylvester matrix corresponding
to the polynomials $\alpha(\xx)$ and $\beta(\xx)$, for elimination of the variable $\basel{x}{n}$, 
is a $(D \times D)$ matrix, where $D = \basel{d}{\alpha}+\basel{d}{\beta}$, with entries either $0$
or any of the multivariate polynomials $\basel{a}{i}(\basel{x}{1}, \, \ldots, \, \basel{x}{n-1})$,
for $0 \leq i \leq \basel{d}{\alpha}$, and $\basel{b}{j}(\basel{x}{1}, \, \ldots, \, \basel{x}{n-1})$,
for $0 \leq j \leq \basel{d}{\beta}$. It is assumed that the number of terms in the sum of terms form
of expansion of the entries is at most $\basel{L}{\max}$ each for these polynomials. The resultant,
denoted by $\mathsf{Res}\bglb\alpha(\xx), \beta(\xx)\bgrb$, with respect to the variable $\basel{x}{n}$,
is the determinant of the $D \times D$ Sylvester matrix. The expansion of the determinant form
as the sum of $D!$ many product terms shows that the number of terms in the resultant can be
at most $D! \basel{L^{D}}{\max} < \bglb D \basel{L}{\max} \bgrb^{D}$.

  Now, comparing with the lower bounds for the number of terms obtained in the case of the
  extended Euclidean algorithms as in Section \ref{Sec-EEA-BA-GB}, it may be found that  
   $\bglb D \basel{L}{\max} \bgrb^{D} $
   $ = $
   $ \basel{L^{D \log_{\basel{L}{\min}} (D \basel{L}{\max})}}{\min}$
   $ = $
   $ \basel{L^{\bglb \frac{D \log_{2} (D \basel{L}{\max})}{\log_{2}(\basel{L}{\min})}\bgrb}}{\min}$.
   With $\basel{d}{\alpha} = \basel{d}{\beta} = d$ and $D = 2d$, if 
   $\bglb \frac{ 2d \log_{2} (2d \basel{L}{\max})}{\log_{2}(\basel{L}{\min})}\bgrb$
   is much smaller than $2^{d}$, then, clearly, $\mathsf{Res}\bglb\alpha(\xx), \beta(\xx)\bgrb$
   is much better. The observation holds, even when $\basel{d}{\alpha} \neq \basel{d}{\beta}$.
   
   On the other hand, if the degree of occurrence of a variable $\basel{x}{i}$ is at most
   $\basel{\Delta}{i}$, for a fixed index $i$, where $1 \leq i \leq n-1$,
   in the multivariate polynomial coefficients of the input polynomials, then
   the degree of occurrence of the variable $\basel{x}{i}$ in 
   ${\mathsf{Res}}\bglb \alpha(\xx), \, \beta(\xx)\bgrb$, with respect to
   the variable $\basel{x}{n}$, is at most $D\basel{\Delta}{i}$.
   With $D = 2d$, if $2d\basel{\Delta}{i}$ is much smaller than
   $2^{d}\basel{\delta}{i}$, as found in the previous section,
   then, again clearly, computation of the parametric resultant
   yields superior performance compared to the extended Euclidean algorithm.
   
   However, the straightforward expansion of the determinant form results in exponential time
   complexity, owing to the $D!$ many terms in the sum. Similarly, the Gaussian elimination
   method could deliver a worse performance than the extended Euclidean algorithm, because the
   intermediate results may not collapse into a small number of terms until the final result.
   
   The following propagation of computations of the determinants of smaller dimension square
   submatrices to larger square submatrices is useful. For the computation of the determinant
   of a $D \times D$ matrix, for a large dimension $D > 1$, let the determinant
   $\det\bglb \basel{S}{k} \bgrb$ and inverse $\basel{S^{-1}}{k}$ of a $k \times k$
   submatrix, $\basel{S}{k}$, be found, where $\basel{S}{k}$ is a submatrix of the
   $k \times D$ matrix $\basel{P}{k}$, obtained by collecting the first $k$ rows of
   the $D \times D$ matrix, inductively, for some $k \geq 2$, but $k \leq D-1$.
   Let  $\basel{P}{k+1}$ be the $(k+1) \times D$ matrix obtained by adjoining
   the next row in the $D \times D$ matrix to $\basel{P}{k}$. Assuming
   that the determinant of the given $D \times D$ matrix, which is
   ${\mathsf{Res}}\bglb \alpha(\xx), \, \beta(\xx)\bgrb$, for the
   Sylvester matrix, with respect to the variable $\basel{x}{n}$,
   does not identically vanish, as an element in the integral domain
   $\polynomials{\scalars}{x}{n-1}$,  the rows of the matrix $\basel{P}{k+1}$
   are linearly independent over the field of  fractions of the integral domain
   $\polynomials{\scalars}{x}{n-1}$. By the equality  of the row rank to the
   column rank, there are $k+1$ linearly independent columns of $\basel{P}{k+1}$.
   Now, of these linearly independent columns, $k$ of the columns can be chosen
   to be those corresponding to the columns of $\basel{S}{k}$, for the following
   reason:   the columns corresponding to $\basel{S}{k}$ are linearly independent,
   by its invertibility, and if every column of $\basel{P}{k+1}$ were a linear combination
   of the $k$ columns, corresponding to those of $\basel{S}{k}$, then the column rank
   of $\basel{P}{k+1}$ itself would be $k$. For the Sylvester matrix with respect
   to the variable $\basel{x}{n}$, the linear combination is taken over the field
   of fractions of the integral domain $\polynomials{\scalars}{x}{n-1}$. Thus,
   at any point, if it is not possible to propagate the computation of the determinant
   from a $k\times k$ submatrix to $(k+1)\times (k+1)$ submatrix, for the reason that
   the column rank cannot increase, after adjoining any of the remaining $D-k$ rows to
   $\basel{P}{k}$, then the determinant of the given matrix itself vanishes, and,
   for the Sylvester matrix, with respect to the variable $\basel{x}{n}$, 
   $\mathsf{Res}\bglb \alpha(\xx), \, \beta(\xx)\bgrb$ itself identically vanishes.
   Given $\basel{S}{k}$, $\det\bglb \basel{S}{k}\bgrb$ and $\basel{S^{-1}}{k}$,
   the computations required for identifying an appropriate column in $\basel{P}{k+1}$,
   in order to form the $(k+1) \times (k+1)$ matrix $\basel{S}{k+1}$, and the determinant and
   the inverse of $\basel{S}{k+1}$, can be performed using standard formulas from 
   matrix algebra. This method of computation of the parametric resultant,
   $\mathsf{Res}\bglb \alpha(\xx), \, \beta(\xx)\bgrb$, avoids needless space
   and time explosion, that can be observed in the Gaussian elimination method.
   
   If $\mathsf{Res}\bglb \alpha(\xx), \, \beta(\xx)\bgrb = 0$, for some interpretation
   of the variables $\basel{x}{i} = \basel{\xi}{i}$ in the algebraic closure of $\scalars$,
   for $1 \leq i \leq n-1$, and neither of $\alpha(\xx)$ and $\beta(\xx)$ vanishes identically
   as respective single variable polynomials in $\basel{x}{n}$, for the ground instances of
   $\basel{x}{i} = \basel{\xi}{i}$, for $1 \leq i \leq n-1$, then $\alpha(\xx)$ and $\beta(\xx)$
   share a common zero in the algebraic closure of $\scalars$. This property is called the
   {\em{minimality}} of the reduction step, for the elimination of the variable $\basel{x}{n}$,
   from the two  participating multivariate polynomials. The collection of multivariate
   polynomials as obtained by computing the parametric resultant of two multivariate
   polynomials, at a time, with respect to any of the independent variables, is called
   a Rabin basis.

{ \textcolor{titlecolor}{\textbf{\section{\label{Sec-PSPACE-not-equal-to-NP-or-P}Proof Showing That $\PSPACE$ Can Be Neither $\NPTIME$ Nor $\PTIME$}}} }
{ \textcolor{titlecolor}{\textbf{\subsection{\label{SubSec-NSPACE-With-Proof}Proof Showing that $\NSPACEwithProof = \NPTIME$}}} }
For the definitions of the computational complexity classes denoted by
$\PTIME$, $\NPTIME$, $\NSPACE$ and $\PSPACE$, the readers are referred
to \cite{HMU-2007}, where $\NSPACE$ appears as ${\mathcal {NPS}}$ and $\PSPACE$
as $\mathcal{PS}$. As an additional complexity class, let $\NSPACEwithProof$
be the collection of languages over an alphabet containing at least two symbols,
such that for the acceptability of an input word, for each language independently,
a nondeterministic requiring space bounded by a polynomial (specific to the
particular language) in the string length of the input word, as with 
$\NSPACE$, exists, but with an additional property that a proof of
acceptance can be automatically generated with respect to any particular
system of logic with its own rules of inference. Some of the rules of
inference may be specialized to the particular language. The proof is assumed
to be machine checkable, for its validity, in $\PSPACE$.  The time required to
check the validity of each step in the proof is taken to be bounded by a fixed,
but sufficiently large, constant. The syntax checking of the proof is also assumed
to require time bounded by a polynomial in the size of the proof. More specifically,
the decision problem of proof checking is required to belong to the complexity
class $\PTIME$. It is easy to see that $\NPTIME \subseteq \NSPACEwithProof$,
because of the deterministic polynomial time verification condition for the
languages in $\NPTIME$.

For the converse inclusion, let $\Sigma$ be an alphabet of at least two distinct symbols,
and let ${\mathcal{L}} \subseteq \Sigma^{\star}$ be a language in $\NSPACEwithProof$.
By assumption, there is a polynomial $\basel{p}{\mathcal{L}}(|\omega|)$, for every word
$\omega \in \Sigma^{\star}$, where $|\omega|$ denotes the string length of $\omega$,
such that whenever $\omega \in {\mathcal L}$, there is a proof attesting to this fact,
of at most $\basel{p}{\mathcal{L}}(|\omega|)$ bits of information, relative to a particular
fixed system of logic, together with the rules of inference, perhaps specialized for the
the language ${\mathcal{L}}$. Now, for every $\omega \in {\mathcal{L}}$, the proof
that $\omega$ indeed belongs to ${\mathcal L}$ can be guessed, checked for syntactic
correctness of the proof, and finally checked for the validity of the proof itself,
in overall time bounded by some polynomial in $|\omega|$. For an external user, 
$\basel{p}{\mathcal{L}}(|\omega|)$ may still remain oblivious, as the technicality
of asserting the membership of ${\mathcal L}$ to $\NSPACEwithProof$ assumes only
its existence.

{ \textcolor{titlecolor}{\textbf{\subsection{\label{SubSec-A-Specific-Computational-Problem}A Specific Computational Problem in $\PSPACE$ Belonging to Neither $\NPTIME$ Nor $\PTIME$}}} }
Let $\fieldChar$ be a large prime number, and $\finiteIntegerField$ be the finite field of
integers with arithmetic operations $\mod \fieldChar$. Let $m, \, n \geq 5$ be positive integers
and $f(t,\, x) = \sum_{i = 0}^{n-1} \basel{a}{i}(t) x^{i} + x^{n} \in 
\singlevariablepolynomials{\finiteIntegerField}{t,\, x}$, where
$\basel{a}{i}(t) \in \singlevariablepolynomials{\finiteIntegerField}{t}$,
for $0 \leq i \leq n-1$, are nonzero polynomials, with undermined coefficients
for the purpose of description of the computational problem, of degree at most
$m$ each. For convenience, it is assumed that there is exactly one polynomial
of degree $m$ and that all the remaining $n-1$ polynomials are of degree at most $m-1$.

The computational question is, ``what is the number of distinct values of
$t$ in the algebraic closure of $\finiteIntegerField$, for which there can
be a solution, for $x \in \finiteIntegerField$, such that $f(t,\, x) = 0$?''
This question is akin to the problems studied along the lines of \cite{Weil-1949}.
Since the polynomials $\basel{a}{i}(t)$, for $0 \leq i \leq n-1$,
are all of degree at most $m$, each, an immediate observation is
to search for values of $t$ in extension fields of $\finiteIntegerField$
of degree at most $m$. As a refinement based on this observation, the computational
task is to enumerate (produce as output) the number of solutions for $t$ in
the extension field of degree $d$ over $\finiteIntegerField$, for each
degree $d$, where $1 \leq d \leq m$, such that $f(t, \, x) = 0$, for some
$x \in \finiteIntegerField$. Obviously, the computational problem is in
$\PSPACE$. In \cite{CGHMP-2003}, the authors discuss another similar
computational problem, but do not assert whether the problem they study
indeed belongs to $\PSPACE$. Instead, their contention is restricted to
the hardness of solving simultaneous multivariate polynomial equations,
in general.

 It is easily observed that $g(t) = {\mathsf {Res}} \bglb f(t,\, x), \, x^{\fieldChar}-x)$
 is a polynomial of degree exactly $m\fieldChar$, and does not depend on $n$,
 by the convenient assumption made.  Also, if $g(t)= 0$, nether of $f(t,\, x)$
 and $(x^{\fieldChar}-x)$ vanishes,  since the leading coefficients of $x$
 in both the polynomials are equal to the constant $1$. However, the multiplicity
 of occurrence of a root of $g(t)$ must be accurately accounted for, as may be
 inferred by observing that the two polynomials $\bglb \phi(t) + x \bgrb$ and
 $\bglb \phi(t) +  (\psi(t))^{2} + x \bgrb$ share a common zero in the
 algebraic closure of $\finiteIntegerField$, exactly when $(\psi(t))^{2} = 0$,
 but each such value of $t$ in the algebraic closure of $\finiteIntegerField$
 must be taken into account as occurring with only multiplicity one.
 Thus, $h(t) = \gcd\bglb g(t), \, g'(t)\bgrb$, where $g'(t)$ is the
 formal derivative of the polynomial $g(t)$, must be computed, and
 finally, the degree of $t$ in the polynomial $\frac{g(t)}{h(t)}$
 yields the answer to the first question. As to the second question,
 the degree of $\gcd$ of  $\frac{g(t)}{h(t)}$ with $(t^{\fieldChar^{d}}-t)$
 yields the number of values of $t$ in the degree $d$ extension of
 $\finiteIntegerField$, such that $f(t, \, x) = 0$, for some
 $x \in \finiteIntegerField$, for $1 \leq d \leq m$.

  It is quite a simple matter to generalize the problem to higher dimensions.
  Let $f(\xx) \in \polynomials{\finiteIntegerField}{x}{n}$, for some integer
  $n \geq 3$, but requiring the number of zeros in the extensions of degree
  at most $\basel{m}{i}$ larger than $2$, for $1 \leq i \leq n$, possibly
  including zeros in $\prod_{i = 1}^{n} GF(\fieldChar, \, \basel{d}{i})$,
  for all possible index vectors $(\basel{d}{1}, \, \ldots, \, \basel{d}{n})$,
  where $1 \leq \basel{d}{i} \leq \basel{m}{i}$ and $GF(\fieldChar, \, \basel{d}{i})$
  is the $\basel{d}{i}$ degree extension of $\finiteIntegerField$. It is worth noting
  that a comparison of $\basel{m}{i}$ to the degree of occurrence of $\basel{x}{i}$
  is omitted, for the purpose of stating the problem in its most generality. The
  enumeration problem can be easily shown to be in $\PSPACE$, because $\basel{m}{i}$,
  for $1 \leq i \leq n$, are fixed inputs to the instance. If it is required to consider
  values for $\basel{x}{i}$ in degree $\basel{d}{i}$ extensions of $\finiteIntegerField$,
  that are not in any smaller dimension extension, then the condition as to whether or not
  $\basel{x^{\fieldChar^{j}}}{i}-\basel{x}{i} \neq 0$ holds, for $1 \leq j \leq \basel{d}{i}-1$,
  must be checked for.  Now, in the algebraic proof, if such a condition must be expressed,
  then the principle of inclusion-and-exclusion must also be applied, in addition to the
  division by $\gcd$ with derivatives, as may be required, for example, in the expression
  for the M\"{o}bius inversion formula. An algebraic proof of validity of the output of the
  $\PSPACE$ algorithm for this problem --- {\em i.e.}, proof of certification for the validity
  of its output as a specific indicator of an algebraic expression --- would definitely
  need an enormously large space.
  
    By accommodating more variables and equations and raising the same question
    concerning the number of distinct solutions to the systems of simultaneous
    multivariate equations,  in general, the contention that $\PSPACE \neq \PTIME$
    can be more  aptly testified, because the solution space cannot be bounded by
    a polynomial, disallowing any claim of producing a direct polynomial time
    algebraic proof, for the validation of the answer produced by a nondeterministic
    algorithm.

  The philosophic question under investigation is whether there can be a shorter form
  of solution for the computational problem exhibited to be in $\PSPACE$. The answer
  to the question is that the algebraic form, depicted as the solution to the problem,
  holds for all the prime numbers occurring in place of the field characteristics and
  for all the other indeterminate parameters as part of the problem instances, and
  hence there cannot be a shorter form, by Herbrand's theorem. It may be recalled
  that, in the first illustration with two variables $t$ and $x$, the degree of
  the resultant of $f(t, \, x)$  and $(x^{\fieldChar}-x)$ does not depend on the
  degree $n$ of occurrence of the variable $x$ in the polynomial $f(t, \, x)$,
  but there are implicitly infinitely many possible choices for the coefficients of
  the polynomials $\basel{a}{i}(t)$, for $0 \leq i \leq n-1$, for various choices of
  $n$ and $\fieldChar$, with the numerical value of $n$ bounded by some polynomial
  in $\log(\fieldChar)$. This situation should not be confused with the way the
  determinant is computed.  In the case of the determinant, there was an easier
  way to compute it, and in the context of the specific computational problem
  in $\PSPACE$, there is simply no alternative algebraic proof attesting the 
  validity of the solution produced by the algorithm requiring space bounded by a 
  polynomial in appropriate values of the parameters of the problem instances.

{ \textcolor{titlecolor}{\textbf{\subsection{\label{SubSec-A-Specific-Computational-Problem}A Specific Computational Problem in $\CoNPTIME$ Belonging to Neither $\NPTIME$ Nor $\PTIME$}}} }
Let $\fieldChar$ be a large prime number, and $\finiteIntegerField$ be the finite field of
integers with arithmetic operations $\mod \fieldChar$. Let $m, \, n \geq 5$ be positive integers
and $f(t,\, x) = \sum_{i = 0}^{n-1} \basel{a}{i}(t) x^{i} + x^{n} \in 
\singlevariablepolynomials{\finiteIntegerField}{t,\, x}$, where
$\basel{a}{i}(t) \in \singlevariablepolynomials{\finiteIntegerField}{t}$,
for $0 \leq i \leq n-1$, are undetermined coefficient polynomials in $t$,
of degree at most $m$ each.

The computational decision question is whether or not the following holds:
$f(t,\, x) \neq 0$, for all $t, \, x \in \finiteIntegerField$?
If $f(t, \, x)= 0$, for some $t,\, x \in \finiteIntegerField$,
then the problem is in $\NPTIME$, and if $f(t,\, x) \neq 0$,
for every $t,\, x \in \finiteIntegerField$, then the problem
is in $\CoNPTIME$. The condition that $f(t,\, x) \neq 0$,
for every $t,\, x\in \finiteIntegerField$, is equivalent
to the condition that $\gcd$ of the three polynomials
$f(t, \, x)$, $(x^{\fieldChar}-x)$ and 
$(t^{\fieldChar}-t) $ is $1$.
A proof of certification could be a derivation that
the $\gcd$ is indeed $1$.

 For a single variable polynomial
$\phi(x) \in \singlevariablepolynomials{\finiteIntegerField}{x}$,
by repeated squaring method,  $x^{2^{i}} \mod \phi(x)$,
for  $i = 1,\, 2,\, 3,\, ...,\, \lfloor\log_{2}(\fieldChar)\rfloor$,
 can be computed, from which
$(x^{\fieldChar}-x) \mod \phi(x)$ can also be computed,
and eventually, $\gcd$ of $(x^{\fieldChar}-x) \mod \phi(x)$ 
and $\phi(x)$ can be shown in a derivation sequence.

However, the same algorithm doesn't work with
two variables, and the following system of
simultaneous multivariate polynomial equations
must be considered:
\begin{equation}
f(t,\, , x ) = 0\,, \tab  (x^{\fieldChar}-x) = 0 \, \tab \textrm{and} \tab 
(t^{\fieldChar}-t) = 0 \label{eqn-in-subsec-a-specific-computational-problem}
\end{equation}
Thus, the resultant of $f(t,\, , x)$ with respect to any of the remaining
two polynomials must be computed, and the $\gcd$ of the resultant and the
other equation left out must be computed. This two-step method is
necessitated by the fact that the last two polynomial equations,
{\em{viz}}, $(x^{\fieldChar}-x) = 0$ and $(t^{\fieldChar}-t) = 0$, 
have no variable in common.  The space requirement for the
resultant in the first step cannot be bounded by a polynomial
in $\log(\fieldChar)$.

Thus, in order to testify the condition that  
$ f(t, \, x)$, $(x^{\fieldChar}-x)$ and
$(t^{\fieldChar}-t)$ have no zeros in common,
simultaneously satisfying them, the size of
the proof cannot be bounded in space by some
polynomial in $\log(\fieldChar)$,
and hence $\CoNPTIME \neq \PTIME$.

    By accommodating more variables and equations and raising the same question concerning
    the nonexistence of solutions to the systems of simultaneous multivariate equations,
    in general, the contention that $\CoNPTIME \neq \PTIME$ can be more aptly
    testified, because the solution space cannot be bounded by a polynomial in the
    acceptable parameter values,  and in particular, by a polynomial in $\log(\fieldChar)$,
    in the preceding example, disallowing any claim of producing a direct polynomial time
    algebraic proof, for the validation of the answer produced by any algorithm.
    
Now, even if a deterministic algorithm for the problem
is equipped with nondeterministic choices of size
bounded by some polynomial in $\log(\fieldChar)$,
the size of the proof remains exponential,
and hence $\CoNPTIME \neq \NPTIME$.
The supporting evidence for the problem
instance to be not in $\NPTIME$ is that
there is no deterministic proof of certification,
bounded in size by some polynomial in the 
parameter values, even if
an assumed deterministic algorithm
for it is equipped with the capability
to make nondeterministic choices of
size bounded by some polynomial
in the acceptable values of the 
parameters of the problem instance.

Some of the implications of the assertions
that $\NPTIME \neq \PTIME$ and that
$\CoNPTIME \neq \NPTIME$ are that
$\PSPACE \neq \PTIME$, 
$\PSPACE \neq \NPTIME$ and
$\PSPACE \neq \CoNPTIME$,
and that the polynomial time hierarchy
does not collapse.
    
{ \textcolor{titlecolor}{\textbf{\subsection{\label{SubSec-Univariate-Polynomial-Modular-Arithmetic-Computations}Computational Complexity of Univariate Polynomial Modular Arithmetic Operations for Polynomials with Sparse Nonzero Coefficients}}} }
 Let $\fieldChar$ be a large prime number, $m, \, n \geq 5$ be positive integers,
 $\basel{a}{k-1},\, \basel{i}{k},$
 $ \basel{b}{l-1},\, \basel{j}{l}$
$ \in $
$ \finiteIntegerField \backslash \{0\}$, for $1 \leq k \leq m$ and $1 \leq l \leq n$.
 Let  $f(x) = \basel{a}{0} + \sum_{k = 1}^{m-1} \basel{a}{k} x^{\basel{i}{k}}+ x^{\basel{i}{m}}$
 and   $g(x) = \basel{b}{0} + \sum_{l = 1}^{n-1} \basel{b}{k} x^{\basel{j}{l}}+ x^{\basel{j}{n}}$.
  The computational decision question is whether
  $\bglb \bglb f(x)\bgrb^{\fieldChar-1}-1\bgrb \cdot $
  $\bglb \bglb g(x)\bgrb^{\fieldChar-1}-1\bgrb = 0 $,
  for all $x \in \finiteIntegerField$. 
  The decision problem can also be
  posed as the question that  
   $\gcd\bglb f(x),\, g(x)\bgrb \neq 0$,
   for all $x \in \finiteIntegerField$,
   and belongs to the complexity class 
   $\CoNPTIME$.
   
   If $g(x)$ is taken to be the polynomial
   $\bglb x^{\fieldChar-1}-1 \bgrb$ and 
   $\basel{a}{0} \neq 0$, then the
   decision problem being discussed reduces
   to the question as to whether $f(x) \neq 0$,
   for all $x \in \finiteIntegerField$, and
   becomes equivalent to the question as to
   whether $\gcd \bglb f(x),\, (x^{\fieldChar-1}-1)\bgrb =1$.
   
   Since the degree $\basel{i}{m}$ of $f(x)$ may be
   chosen to be arbitrarily large, the computation of
   $(x^{\fieldChar-1}-1) \mod f(x)$ could require space
   that is exponential in the string length of $f(x)$,
   due to its sparseness. In the derivation sequence of
   $\gcd \bglb f(x),\, (x^{\fieldChar-1}-1)\bgrb$
   the coefficients of intermediate polynomials in
   the proof script can be interpreted as hash
   functions in the undetermined coefficients
   $\basel{a}{k}$, for $0 \leq k \leq m-1$.
   If each of the hash functions, treated independently,
   vanishes with only exponentially insignificant
   probability of only once in $\fieldChar$ many chances,
   for any instantiation of the undetermined coefficients
   $\basel{a}{k}$, for $0 \leq k \leq m-1$, then the
   space requirement for the derivation sequence cannot
   be upper bounded by any polynomial in the string length
   of the sparse polynomial $f(x)$.  This evidence fortifies
   the contention that $\CoNPTIME \neq \PTIME$.
   
   The descriptions of in the next couple of paragraphs
   are intended to illustrate how to navigate possible
   shortcomings. The first example illustrates an attempt
   to generate an instance of $\CoNPTIME$, but ends up
   with an instance in $\NPTIME$, thereby belittling any
   attempt to claim that there cannot be a shorter proof.
   In the second example, a more general problem is
   considered, but the complementary instance cannot be
   guaranteed to belong to $\NPTIME$, thereby restricting
   the application of its scope to the contention that
   $\PSPACE \neq \PTIME$. The first example is rarefied
   subsequently, in a third paragraph, illustrating how
   to indirectly show the nonexistence of a shorter proof.
   
    In order to show the existence of sparse
   polynomials  in $\singlevariablepolynomials{\finiteIntegerField}{x}$,
   that do not vanish for any interpretation of $x \in \finiteIntegerField$,
   let $\basel{\gamma}{i}$,  $\basel{\delta}{i}$,  $\basel{r}{i}$
   $ \in \finiteIntegerField\backslash \{0\}$
   such that $\gcd(\basel{r}{i},\, p-1) >  1$
   and $\bglb \basel{\gamma^{-1}}{i}\basel{\delta}{i} \bgrb^{\frac{\fieldChar-1}{\gcd(\basel{r}{i},\, p-1)}} \neq 1 \mod \fieldChar$,
   for $1 \leq i \leq \kappa$, for some fixed small integer $\kappa$, such as about $20$.
   Some of the integers $\basel{r}{i}$ can be conveniently chosen to be larger than
   $\frac{\fieldChar}{2}$, where $1 \leq i \leq \kappa$.
   The polynomial $\prod_{i = 1}^{\kappa} \bglb \basel{\gamma}{i} x^{\basel{r}{i}} - \basel{\delta}{i} \bgrb$,
   after expansion and simplification applying the rule $x^{\fieldChar-1} \equiv 1 \mod \bglb x^{\fieldChar-1}-1\bgrb$,
   is taken to be a polynomial $f(x)$ in the simplified form.
   It is convenient to choose $\basel{r}{i}$ such that there are ample number of identities
   of the form $\basel{r}{i}+\basel{r}{j} = \basel{r}{k}+\basel{r}{l} \mod (\fieldChar-1)$,
   for some distinct indexes $i, \, j, \, k, \, l$, and more generally, the indexes may be
   carefully selected with an objective of keeping the number of terms in the polynomial
   $f(x)$ small.  However, the instances of $f(x)$ as constructed may also be shown to belong
   to the complexity class $\NPTIME$, as well, by observing that the hidden parameters may be
   guessed and substituted in the corresponding expressions, to reconstruct the polynomial $f(x)$.
   In this case, there is a possibility of claiming the existence of a shorter proof to the
   effect of showing that $f(x) \neq 0$, for every $x \in \finiteIntegerField$, by checking
   the derivation of $f(x)$ as assumed.  The attempts in the preceding subsections are intended
   to  somehow indirectly show that this hope of generating an easy proof does not exist.
   
   One more passing insight is that if the question in this section was to show
   $\gcd\bglb f(x),\, g(x)\bgrb = 1$, then it is not easy to show that the problem
   belongs to $\CoNPTIME$, because the complementary problem of checking whether
   $f(x)$ and $g(x)$ have a nontrivial common factor, or equivalently, whether
   they share a common zero in the algebraic closure of $\finiteIntegerField$,
   is not guaranteed to be in $\NPTIME$. The irreducible common factors of
   $f(x)$ and $g(x)$ may be of large degree, such as comparable to some constant
   fraction of $\fieldChar$. Thus, in this case, the contention is restricted
   to $\PSPACE \neq \PTIME$.
   
   The first example is rarefied in this paragraph. Let $f(x)$ be
   a sparse polynomial of arbitrary large degree, satisfying the
   Eisenstien's criterion \cite{Lang-2002}, for irreducibility
   over the field of rational numbers. Constructively, let
   $f(x) = \basel{a}{0} + \sum_{k = 1}^{m-1} \basel{a}{k} x^{\basel{i}{k}} + x^{\basel{i}{m}}$,
   such that $1 \leq \basel{i}{1} < \ldots < \basel{i}{m} < \fieldChar-1$,
   where the coefficients $\basel{a}{k}$, $0 \leq k \leq m-1$
   are chosen as follows: let $\primeInt$ be a large prime integer,
   different from $\fieldChar$, $\basel{a'}{k}$ be arbitrary, 
   but large integers, such that  $\basel{a}{k} = \primeInt \basel{a'}{k}$
   and $\gcd \bglb \basel{a}{0},\, \primeInt\bgrb$
   $ = $
   $\gcd \bglb \basel{a'}{k},\, \primeInt\bgrb$
   $ = $
   $1 \mod \primeInt$,   for $1 \leq k \leq m-1$.
   The polynomial $f(x)$ is identified with $f(x) \mod \fieldChar$,
   by simplifying its coefficients in $\finiteIntegerField$.
   Now, it is tactically assumed that, with some positive
   probability,  $\gcd\bglb f(x), (x^{\fieldChar}-x) \bgrb = 1$,
    If this condition must be deterministically checked,
    during the construction phase, then the route followed
    is as follows:  a small degree sparse polynomial $h(x)$ is
    constructed, such that $\gcd\bglb h(x),\, (x^{\fieldChar}-x)\bgrb =1$,
    and $f(x)$ is taken to be the polynomial obtained by
    substituting $x^{r}$ for $t$ in $h(t)$, where $r$
    is a large positive integer relatively prime to $(\fieldChar-1)$.
    The exponents in $f(x)$ are taken $\mod (\fieldChar-1)$.
    As discussed before, the possibility of 
    deterministic polynomial time proof checking,
    after allowing for nondeterministic choices
    of hidden nature, relinquishes any claim
    that there cannot be short deterministic proof.
    With the tactical assumption, $f(x)$ is only
    probabilistically assumed to be relatively prime
    to $(x^{\fieldChar}-x)$, but there is no direct
    evidence of a short deterministic proof
    testifying the claim that $\gcd$ of $f(x)$
    and  $(x^{\fieldChar}-x)$ is indeed $1$.
    The tactical assumption is contented with
    some positive probability, however little
    and insignificant it may be, as may be seen
    to be different from the restriction posed
    for a decision problem to be in the complexity
    class $\NPTIME$, {\em{viz}}, that the verification
    step must comply with deterministic polynomial time 
    proof checking. Even if the conditions of the
    Eisenstein's criterion are deterministically shown
    to be satisfied, it may be observed that,  in the
    construction phase, any attempt to deterministically
    check that $\gcd\bglb f(x), (x^{\fieldChar}-x) \bgrb = 1$
    was deliberately avoided, for precisely this reason.
    As a last refinement, it is still possible to pose
    the question as to whether 
    $\gcd\bglb f(x), (x^{\nu}-1) \bgrb = 1$,
    where $\nu$ is a large positive integer
    dividing $(\fieldChar-1)$. The probability
    of the event $\gcd\bglb f(x), (x^{\nu}-1) \bgrb = 1$,
    for a large divisor $\nu$ of $(\fieldChar-1)$,
    which is strictly smaller than $(\fieldChar-1)$ itself,
    is enhanced considerably when compared to the
    probability of the event 
     $\gcd\bglb f(x), (x^{(\fieldChar-1)}-1) \bgrb = 1$.
    This prospect of enhanced probability effectively 
    seals the contention that $\CoNPTIME \neq \PTIME$.

{ \textcolor{titlecolor}{\textbf{\subsection{\label{SubSec-IP}Implications of the Fact that $\PSPACE = \IP$}}} }
In \cite{Shamir-1992}, it is shown that $\IP = \PSPACE$.
By Savitch's theorem \cite{HMU-2007}, the complexity classes
 $\NSPACE$ and $\PSPACE$ represent the same complexity class.
 Thus, $\PSPACE$ is closed under augmentation by the capability
 to make nondeterministic choices of size bounded by some
 polynomial in the string length of the input word.
Allowing for the nondeterministic choices of space bounded
by a polynomial in the string length of the input word, 
the class of languages that admit bounded error probabilistic
polynomial time proofs, as computable by nondeterminstic polynomial
time algorithms, is exactly $\IP$. 
Let $\BoundedErrorProbPTIME$ be the class of languages acceptable
by bounded error probabilistic algorithms with probabilistic
polynomial time proofs attesting the membership of an input word
for each such language. Tentatively, if it is assumed that
$\BoundedErrorProbPTIME = \PTIME$, then clearly, it must be
the case that $\IP = \NSPACE = \PSPACE = \NPTIME$.
However, the discussion of the previous subsection
shows that $\PSPACE \neq \NPTIME$, implying that 
$\BoundedErrorProbPTIME \neq \PTIME$. In fact,
by equipping a probabilistic algorithm with the
capability to make nondeterministic choices of
size bounded by some polynomial in the string length
of the input word,  an algorithm for a language in
the complexity class $\NSPACE$ can be emulated,
but a nondeterminstic polynomial time algorithm,
equipped  with the same additional capability,
remains only a nonderministic polynomial time
algorithm. Thus, $\BoundedErrorProbPTIME \neq \NPTIME$,
and by complementarity, $\BoundedErrorProbPTIME \neq \CoNPTIME$,
as well. By the same reasoning, if a probabilistic algorithm
for some language in the complexity class $\BoundedErrorProbPTIME$
is equipped with an additional capability of
a polynomial time algorithm, then it still remains
a probabilistic algorithm for the some or some
other language in the complexity class $\BoundedErrorProbPTIME$.
But if the probabilistic algorithm is equipped with
the capability of nondeterministic choices
bounded by some polynomial in the string
length of the input word, then the said
algorithm effectively emulates an algorithm
for some language in the complexity class $\IP$,
and therefore, clearly, $\NPTIME \neq \PTIME$,
whereby it follows that $\CoNPTIME \neq \PTIME$,
as well. The informal arguments can be made
precise by introducing appropriate
relativization of complexity classes.

{ \textcolor{titlecolor}{\textbf{\section{\label{Sec-Concl}Conclusions}}} }
This paper presents an original reduction method for solving simultaneous multivariate 
polynomial equations, by eliminating one variable, from two equations, taken at a time.
The reduction method is shown to satisfy a certain minimality criterion, and hence
becomes optimal in respect of the constraints stated. A mathematical problem that
is in $\PSPACE$ but that which cannot be in either of $\NPTIME$ and $\PTIME$ is also
presented. As an afterthought, the trace of execution of an algorithm, even allowing
for nondeterministic choices of sizes bounded by some polynomial in the appropriate
values of the input parameters for the instances, must also be bounded in size
by some polynomial in those parameter values. The trace may be supplied as the
input to a debugger program for validation of its operation. This part may be
included in the polynomial time certification for the algorithm, for each
given input instance. The contentments assert that $\PSPACE \neq \NPTIME$
and $\PSPACE \neq \PTIME$, affirmatively. For probabilistic algorithms,
the class $\BoundedErrorProbPTIME$ of languages decidable by bounded error
probabilistic algorithms with probabilistic polynomial time proofs for the
membership of an input word is not the same as the complexity class $\PTIME$.
In fact, by equipping a probabilistic algorithm with the
capability to make nondeterministic choices of size bounded
by some polynomial in the string length of the input word,
an algorithm for some language in $\NSPACE$ can be emulated,
but a nondeterminstic polynomial time algorithm, equipped
with the same additional capability, remains only a nonderministic
polynomial time algorithm. Thus, $\BoundedErrorProbPTIME \neq \NPTIME$, as well.
It follows that $\NPTIME \neq \PTIME$ --- hence that $\CoNPTIME \neq \PTIME$ ---
by relativization of $\BoundedErrorProbPTIME $, with respect to $\PTIME$ and $\NPTIME$.

{  \textcolor{titlecolor}{\textbf{\section*{\label{Sec-ACK}Acknowledgements}}} }
\textcolor{secheader}{
The authors gratefully acknowledge fruitful discussions with Professor Michael Oser Rabin, Professor Ronald Linn Rivest and Professor Adi Shamir. The second author received introductions about machine checkable proofs and polynomial time verification certifications during a stint of postdoctoral position under the supervision of Professor Amir Pnueli, at the Weizmann Institute of Science, during the year 2002.
}
\\
\\
{
\textcolor{refscolor}{

}
}
\end{document}